\UseRawInputEncoding
\documentclass[11pt,twoside, leqno]{article}
\usepackage{mathrsfs}
\usepackage{amssymb}
\usepackage{amsmath}
\usepackage{amsthm}
\usepackage{color}

\usepackage[active]{srcltx}
\usepackage{amsfonts}
\usepackage[notref,notcite]{showkeys}
\allowdisplaybreaks

\def\supp{{\rm{\,supp\,}}}

\def\bint{{\ifinner\rlap{\bf\kern.30em--}
\int\else\rlap{\bf\kern.35em--}\int\fi}\ignorespaces}

\def\sbint{{\ifinner\rlap{\bf\kern.32em--}
\hspace{0.078cm}\int\else\rlap{\bf\kern.45em--}\int\fi}\ignorespaces}

\newtheorem{theorem}{Theorem}[section]
\newtheorem{lemma}[theorem]{Lemma}

\theoremstyle{definition}
\newtheorem{remark}[theorem]{Remark}
\newtheorem{definition}[theorem]{Definition}
\numberwithin{equation}{section}

\numberwithin{equation}{section}

\numberwithin{equation}{section}

\begin{document}

\arraycolsep=1pt

\title{\Large\bf Bilinear $\theta$-type Calder\'on-Zygmund operators and its commutator on generalized weighted Morrey spaces over RD-spaces\footnotetext{\hspace{-0.35cm}
\endgraf Author and E-mail address:\ Suixin He,\ hesuixinmath@126.com
\endgraf*Corresponding author and E-mail address:\ Shuangping Tao,\ taosp@nwnu.edu.cn
}}
\author{Suixin He and Shuangping Tao$^{\ast}$\\[3pt]
\small (College of Mathematics and Statistics, Northwest Normal University, Lanzhou, Gansu, 730070,\\[3pt]\small  P. R. China)}
\date{ }
\maketitle

\begin{center}
\begin{minipage}{13cm}\small
{\noindent{\bf Abstract:}}\ \ An RD-space $\mathcal{X}$ is a space of homogeneous type in the sense of Coifman and Weiss with the additional property that a reverse doubling property holds in $\mathcal{X}$. In this setting, the authors establish the boundedness of bilinear $\theta$-type Calder\'on-Zygmund operator $T_{\theta}$ and its commutator $[b_1,b_2,T_{\theta}]$ generated by the function $b_1,b_2\in BMO(\mu)$ and $T_{\theta}$ on generalized weighted Morrey space $\mathcal{M}^{p,\phi}(\omega)$ and generalized weighted weak Morrey space $W\mathcal{M}^{p,\phi}(\omega)$ over RD-spaces.
\\[2pt]

\mbox{}\textbf{Keywords:}\ RD-space$\cdot$ Bilinear $\theta$-type Calder\'{o}n-Zygmund operator $\cdot$ Commutator $\cdot$ Generalized weighted Morrey spaces \\[2pt]

\mbox{}\textbf{2020 MR Subject Classification:}\ \ Primary 43A85; Secondary 42B20 $\cdot$ 42B35\\
\end{minipage}
\end{center}


\section{Introduction}
\quad\
The space of homogeneous type, first introduced by Coifman and Weiss \cite{CW1,CW2}, is a general framework for studying the Calder\'on-Zygmund operators and functions spaces. Around 1970s, Coifman and Weiss began to investigate the some classical harmonic analysis problems on the metric space, which is called space of homogeneous type $(\mathcal{X},d,\mu)$,  equipped with a metric $d$ and a regular Borel measure $\mu$ satisfying the doubling condition, if there exists a positive constant $C_{\mu}>1$ such that, for any ball $B(x,r):=\{y\in \mathcal{X}:d(x,y)<r\}$ with $x\in \mathcal{X}$ and $r>0$ ,
\begin{eqnarray}
\mu(B(x,2r))\leq C_{\mu}\mu(B(x,r))
\end{eqnarray}
holds. Since then, many  classical results were extended to the spaces of homogeneous type in the sense of Coifman and Weiss. However, some results have so far obtained only on the RD-spaces, which
means that $(\mathcal{X},d,\mu)$ is a space of homogeneous type if there exists positive constants $a,b>1$ such that,
\begin{eqnarray}
b\mu(B(x,r))\leq \mu(B(x,ar));
\end{eqnarray}
holds for all $x\in \mathcal{X}$ and $r\in (0,\rm diam(\mathcal{X})/a)$. On the development and research of the operators over RD-spaces, we refer readers see, e.g., \cite{GLY,JYY,LY,YZ,ZLL}.

In recent years, solving numerous important problems in harmonic analysis and PDEs became possible due to the progress reached in the weighted theory, generally speaking, in new function spaces; see, for instance, \cite{P,KMRS}. Moreover, the $A_{p}$ weight theory, which was first studied in \cite{M1}, is one of the cores of the weighted theory.
It should be pointed out that Morrey space, which were introduced by Morrey in 1938 (see \cite{M}) in order to study regularity questions which appear in the calculus of variations, describe local regularity
more precisely than Lebesgue spaces and widely use not only harmonic analysis but also partial differential equations (see \cite{G,GT}).
In 2009, Komori and Shirai \cite{KS} introduced the weighted Morrey spaces and study the several properties of classical operators on the classical Euclidean space. The further research and development about the Morrey spaces and weighted Morrey spaces over different settings, the readers can see \cite{A,CF,DR,DXY,HZ,N,P1,S,SE,ST} and the references therein. In 2020, Chou et al. \cite {CL} introduced the generalized weighted Morrey spaces over RD-spaces, as applications, the boundedness of the classical operator was established. Very recently, Li et al. \cite{LLW} established that the boundedness of the commutators generalized by the $\theta$-Calder\'{o}n-Zygmund operator and the BMO functions in generalized weighted Morrey spaces over RD-spaces.

Motivated by the above research, in this article, we will mainly study the boundedness of bilinear $\theta$-type Calder\'{o}n-Zygmund operator $T_{\theta}$ and its commutator $[b_1,b_2,T_{\theta}]$ associated with function $b_1,b_2\in BMO(\mu)$ on generalized weighted Morrey space $\mathcal{M}^{p,\phi}(\omega)$ and generalized weighted weak Morrey space $W\mathcal{M}^{p,\phi}(\omega)$ over $(\mathcal{X},d,\mu)$. Before present the mainly results of this paper, we first recall some necessary definitions and notion.

The definition of the generalized weighted Morrey space introduced by Chou et al. \cite{CL}.
\begin{definition}
Let $p\in[1,\infty)$, $\phi:(0,\infty)\rightarrow (0,\infty)$ be an increasing functions and $\omega$ be a weight on $\mathcal{X}$. Then the generalized weighted Morrey space $\mathcal{M}^{p,\phi}(\omega)$, equipped with the norm
$$\|f\|_{\mathcal{M}^{p,\phi}(\omega)}:=\mathop{\sup}\limits_{B}\bigg\{\frac{1}{\phi(\omega(B))}\int_{B}|f(x)|^{p}\omega(x)\mathrm{d}\mu(x)\bigg\}^{\frac{1}{p}}<\infty.$$

We also denote by $W\mathcal{M}^{p,\phi}(\omega)$ the generalized weighted weak Morrey space of all locally integrable functions $f$ satisfying
$$\|f\|_{W\mathcal{M}^{p,\phi}(\omega)}:=\mathop{\sup}\limits_{B}\mathop{\sup}\limits_{t>0}\frac{1}{[\phi(\omega(B))]^{\frac{1}{p}}}t\omega(\{x\in B:|f(x)|>t\})^{\frac{1}{p}}<\infty.$$
\end{definition}

\begin{remark}
(i) When $\phi(x)=1$, $\mathcal{M}^{p,\phi}(\omega)=L^{p}(\omega)$ and $W\mathcal{M}^{p,\phi}(\omega)=WL^{p}(\omega)$. Consequently, the generalized weighted (weak) Morrey space is an extension of the weighted (weak) Lebesgue space.

(ii) When $\phi(x)=x^{k}$ with $0<k<1$, then $\mathcal{M}^{p,\phi}(\omega)=\mathcal{M}^{p,k}(\omega)$ and $W\mathcal{M}^{1,\phi}(\omega)=W\mathcal{M}^{1,k}(\omega)$. Hence, the generalized weighted (weak) Morrey space is an extension of the weighted (weak) Morrey space.
\end{remark}

In what follows, let $V(x,y):=\mu(B(x,d(x,y)))$, now we state the definition of bilinear $\theta$-type Calder\'{o}n-Zygmund operator as follows.
\begin{definition}
Let $\theta$ is be a non-negative nondecreasing functions on $(0,\infty)$ satisfy the following Dini condition:
\begin{eqnarray}
\int^{1}_{0} \frac{\theta(t)}{t}dt<\infty.
\end{eqnarray}

A kernel $K(\cdot,\cdot,\cdot)\in L^{1}_{loc}(\mathcal{(X)}^{3}\backslash\{(x,y_1,y_2):x=y_1=y_2\})$ is called the bilinear $\theta$-type Calder\'{o}n-Zygmund kernel if there exists a positive constant $C$, such that

(1) for all $(x,y_1,y_2)\in \mathcal{X}^{3}$ with $x\neq y_{j}$ for $j\in\{1,2\}$,
\begin{eqnarray}
|K(x,y_1,y_2)|\leq C \frac{1}{[\sum^{2}_{j=1}V(x,y_{j})]^{2}};
\end{eqnarray}

(2) there exists a positive constants $C$, such that, for all $x,x^\prime,y_1,y_2$ with $d(x,x^\prime)\leq c\max_{1\leq j\leq 2}d(x,y_j)$,

\begin{eqnarray}
|K(x,y_1,y_2)-K(x^{\prime},y_1,y_2)|\leq C\theta\bigg(\frac{d(x,x^\prime)}{\sum^{2}_{j=1}d(x,y_j)}\bigg)\frac{1}{[\sum^{2}_{j=1}V(x,y_{j})]^{2}}.
\end{eqnarray}
\begin{remark}\rm
Let $\delta\in (0,1]$, $\theta(t)=t^{\delta}$ with $t>0$, then $K(x,y_1,y_2)$ defined as in Definition 1.3 is just the standard Calder\'{o}n-Zygmund kernel.
\end{remark}
Let $L_b^{\infty}(\mu)$ be the space of all $L^{\infty}(\mu)$ functions with bounded support. An operator $T_{\theta}$ is called a bilinear $\theta$-type Calder\'{o}n-Zygmund operator with $K$ satisfying (1.4) and (1.5) if, for all $f_1,f_2\in L_b^{\infty}(\mu)$ and $x\notin \bigcap^{2}_{j=1} \supp~{f_j}$,

\begin{eqnarray}
T_{\theta}(f_1,f_2)(x)=\int_{\mathcal{X}}\int_{\mathcal{X}}K(x,y_1,y_2)f_{1}(y_1)f_{2}(y_2)d\mu(y_1)d\mu(y_2).
\end{eqnarray}

Given $b_1,b_2\in{BMO}(\mu)$, the commutators $[b_1,b_2,T_{\theta}]$ associated with the bilinear $\theta$-type Calder\'{o}n-Zygmund operator $T_{\theta}$ is respectively defined by
\begin{eqnarray*}
[b_1,b_2,T_{\theta}](f_1,f_2)(x)&=&b_1(x)b_2(x)T_{\theta}(f_1,f_2)(x)-b_1(x)T_{\theta}(f_1,b_2f_2)(x)\\
\quad &&-b_2(x)T_{\theta}(b_1f_1,f_2)(x)+T_{\theta}(b_1f_1,b_2f_2)(x).
\end{eqnarray*}
Also, $[b_1,T_{\theta}]$ and $[b_2,T_{\theta}]$ are defined as follows,
\begin{eqnarray}
[b_1,T_{\theta}](f_1,f_2)(x):=b_1(x)T_{\theta}(f_1,f_2)(x)-T_{\theta}(b_1f_1,f_2)(x),
\end{eqnarray}
\begin{eqnarray}
[b_2,T_{\theta}](f_1,f_2)(x):=b_2(x)T_{\theta}(f_1,f_2)(x)-T_{\theta}(f_1,b_2f_2)(x).
\end{eqnarray}
\end{definition}

Now let's recall the definition of $BMO(\mu)$ space.
\begin{definition}\rm
A function  $b\in L^{1}_{loc}(\mu)$ is said to be in the space $BMO(\mu)$, equipped with the norm
$$\|b\|_{BMO(\mu)}=\mathop{\sup}\limits_{B}\frac{1}{\mu(B)}\int_{B}|b(x)-b_{B}|d\mu(x)<\infty,$$
where the supremum is taken over all balls $B\subset\mathcal{X}$ and
$$b_{B}=\frac{1}{\mu(B)}\int_{B}b(y)\mathrm{d}\mu(y).$$
\end{definition}

\begin{definition}\rm
A weight $\omega$ is said to belong to the class $A_{p}$ for $p>1$ if for any $B\subset\mathcal{X}$,
$$\mathop{\sup}\limits_{B}\bigg(\frac{1}{\mu{(B)}}\int_{B}\omega(y)d\mu(y)\bigg)\bigg(\frac{1}{\mu(B)}\int_{B}\omega^{1-p^{\prime}}(x)d\mu(y)\bigg)^{p-1}<\infty,$$
and $\omega$ belong to the class $A_{1}$, if there is a constant $C$ such that for any $B\subset\mathcal{X}$,
$$\frac{1}{\mu(B)}\int_{B}\omega(y)d\mu(y)\leq C\mathop{\inf}\limits_{y\in B}\omega(x).$$
We denote $A_{\infty}$ class in the natural way by $A_{\infty}=\bigcup_{p>1}A_{p}.$
\end{definition}
To get the boundedness of $[b,T_{\theta}]$ with $b\in BMO(\mu)$ on the generalized weighted Morrey space $\mathcal{M}^{p,\phi}(\omega)$, Wang \cite{W} also suppose that the function $\phi$ in  Definition 1.1 need following condition: there exist a positive constants $k\in [0,1)$ and $C$ such that
\begin{eqnarray}
\frac{\phi(r)}{r^{k}}\leq C\frac{\phi(s)}{s^{k}}~~ for~ all~~ 0<s\leq r<+\infty.
\end{eqnarray}

In addition, to obtain the boundedness of $T_{\theta}$ with $\theta=t^{\delta}$, that is to say $K(x,y)$ is standard Calder\'{o}n-Zygmund kernel, the following two conditions are used by Chou et al. \cite{CL}: there exist two positive constants $C_1$ and $C_2$ such that
\begin{eqnarray}
\frac{\phi(r)}{r}\leq C_1\frac{\phi(s)}{s}~~for~ all~~ 0<s\leq r<+\infty.
\end{eqnarray}
\begin{eqnarray}
\int_{t}^{\infty}\frac{\phi(s)}{s}\frac{ds}{s}\leq C\frac{\phi(t)}{t}~~ for~ all~~t\in(0,\infty).
\end{eqnarray}
\begin{remark}
The condition (1.9) implies the conditions (1.10) and (1.11), see e.g.\cite{CL}.
\end{remark}

The following statements are our main results.
\begin{theorem} \rm  Assume that $T_{\theta}$ is a bilinear $\theta$-type Calder\'{o}n-Zygmund operator with kernel $K$ satisfying (1.4) and (1.5). Let $p_1,p_2\in [1,\infty)$ with $\frac{1}{p}=\frac{1}{p_1}+\frac{1}{p_2}$, ${\omega}\in A_{p}(\mu)$.
Let $\theta:[0,\infty)\rightarrow [0,\infty)$ be a nondecreasing function satisfying (1.3) and $\phi:(0,\infty)\rightarrow(0,\infty)$ be an increasing, continuous function satisfying conditions (1.10) and (1.11), we have the following:\\
(i) when all $p_i>1$, there exists a constant $C$ such that
$$\|T_{\theta}(f_1,f_2)\|_{\mathcal{M}^{p,\phi}(\omega)}\leq C \|f_1\|_{\mathcal{M}^{p_1,\phi}(\omega)}\|f_2\|_{\mathcal{M}^{p_2,\phi}(\omega)},$$
(ii) when some $p_i=1$, there exists a constant $C$ such that
$$\|T_{\theta}(f_1,f_2)\|_{W\mathcal{M}^{p,\phi}(\omega)}\leq C \|f_1\|_{\mathcal{M}^{p_1,\phi}(\omega)}\|f_2\|_{\mathcal{M}^{p_2,\phi}(\omega)}.$$
\end{theorem}

\begin{theorem} \rm  Let $p_1,p_2\in [1,\infty)$ with $\frac{1}{p}=\frac{1}{p_1}+\frac{1}{p_2}$, ${\omega}\in A_{p}(\mu)$ and $b_1,b_2\in BMO(\mu)$.
Let $\theta:[0,\infty)\rightarrow [0,\infty)$ be a nondecreasing function satisfying (1.3) and $\phi:(0,\infty)\rightarrow(0,\infty)$ be an increasing, continuous function satisfying conditions (1.10) and (1.11), we have the following:\\
(i) when all $p_j>1$, there exists a constant $C$ such that
$$\|[b_1,b_2,T_{\theta}](f_1,f_2)\|_{\mathcal{M}^{p,\phi}(\omega)}\leq C \prod^{2}_{j=1}\|b_j\|_{BMO(\mu)}\|f_j\|_{\mathcal{M}^{p_j,\phi}(\omega)},$$
(ii) when some $p_j=1$, there exists a constant $C$ such that
$$\|[b_1,b_2,T_{\theta}](f_1,f_2)\|_{W\mathcal{M}^{p,\phi}(\omega)}\leq C \prod^{2}_{j=1}\|b_j\|_{BMO(\mu)}\|f_j\|_{\mathcal{M}^{p_j,\phi}(\omega)}.$$
\end{theorem}

Finally, we make some conventions on notation. Throughout the paper, $C$ represents a positive constant being independent of the main parameters involved, but may vary from line to line. For a $\mu$-measurable set $E$, $\chi_{E}$ denotes its characteristic function. The symbol $f\lesssim g$ means that there exists a positive constant $C$ such that $f\leq Cg$. Given a ball $B\subseteq\mathcal{X}$ and $\lambda>0$, $\lambda B$ denote the ball which has the same center of $B$ and the radius is $t$ times of $B$. For any exponent $p>1$, we denote by $p'$ its conjugate index, i.e., $\frac{1}{p}+\frac{1}{p'}=1$.


\section{Preliminaries}
\quad To prove the main theorems of this paper, in this section, we give some auxiliary lemmas.
\begin{lemma}\rm \cite{ST1}
Let $p\in(1,\infty)$ and $\omega\in A_{p}(\mu)$. There exist positive constants $C_1$ and $C_{2}$ such that for any ball $B\subset \mathcal{X}$ and each measurable set $E\subseteq B$,
\begin{eqnarray*}
\frac{\omega(E)}{\omega(B)}\leq C_1 \bigg[\frac{\mu(E)}{\mu(B)}\bigg]^{\frac{1}{p}} \quad and  \quad \frac{\omega(E)}{\omega(B)}\geq C_2 \bigg[\frac{\mu(E)}{\mu(B)}\bigg]^{p}.
\end{eqnarray*}
\end{lemma}

\begin{lemma}\rm \cite{CL}
Let $(\mathcal{X},d,\mu)$ be an RD-space, if $\omega\in A_{p}(\mu)$, $p\in(1,\infty)$, then there exist positive constants $C_3, C_{4}>1$ such that for any ball $B\subset \mathcal{X}$,
\begin{eqnarray}
\omega(2B)\geq C_3  \omega(B),
\end{eqnarray}
\begin{eqnarray}
\omega(2B)\leq C_4  \omega(B).
\end{eqnarray}
\end{lemma}

\begin{lemma}\rm\cite{LLW}
Let $p\in(1,\infty)$, $\omega\in A_{p}(\mu)$ and $\phi: (0,\infty)\rightarrow(0,\infty)$ be an increasing, continuous function satisfying (1.10) and (1.11). Then there exists a positive constant $C$ such that for any ball $B\subset X$,
$$
\sum^{\infty}_{k=1}k\bigg[\frac{\phi(\omega(2^{k}B))}{\omega(2^{k}B)}\bigg]^{\frac{1}{p}}\leq C \bigg[\frac{\phi(\omega(B))}{\omega(B)}\bigg]^{\frac{1}{p}}.$$
\end{lemma}

\begin{lemma}\rm \cite{LLW}
Let $b\in BMO(\mu)$ and $\omega\in A_{p}(\mu)$ with $p\in[1,\infty)$, then
\begin{enumerate}
\item[(1)]  there exists a positive constant $C$ such that, for any ball $B\subset\mathcal{X}$ and $k\in \mathbb{Z}^{+}$,
$$|b_{2^{k+1}B}-b_{B}|\leq C(k+1)\|b\|_{BMO(\mu)};$$
\item[(2)]  there exists a positive constant $C$ such that, for any ball $B\subset\mathcal{X}$
$$\bigg\{\int_{B}|b(x)-b_{B}|^{p}\omega(x)d\mu(x)\bigg\}^{\frac{1}{p}}\leq C\|b\|_{BMO(\mu)}[\omega(B)]^{\frac{1}{p}}.$$
\end{enumerate}
\end{lemma}

Finally, we need to recall the following boundedness of the bilinear $\theta$-type Calder\'{o}n-Zygmund $T_{\theta}$ and commutator $[b_1,b_2,T_{\theta}]$ on the weighted Lebesgue space $L^{p}(\omega)$.
\begin{lemma}
Let $K(\cdot,\cdot,\cdot)$ satisfy (1.4) and (1.5), ${\omega}\in A_{p}(\mu)$ and $\frac{1}{p}=\frac{1}{p_1}+\frac{1}{p_2}$ with $1<p_1,p_2<\infty$. Then $T_{\theta}$ can be extended to a bounded operator from $L^{p_1}(\omega)\times L^{p_2}({\omega})$ to $L^{p}({\omega})$.
\end{lemma}
{\it Proof.} From \cite{GLMY}, it is not hard to get that the $T_{\theta}$ is bounded from weighted Lebesgue space $L^{p_1}(\omega)\times L^{p_2}({\omega})$ to $L^{p}({\omega})$. Hence, here we omit the proof.

By an argument analogous to the proof of \cite{W} with a slight modification, we obtain the following results, for briefly, we omit the details here.
\begin{lemma}
Let $b_1,b_2 \in BMO(\mu)$, ${\omega}\in A_{p}(\mu)$. Then commutators $[b_1,b_2,T_{\theta}]$ is bounded from the product of space $L^{p_1}(\omega)\times L^{p_2}({\omega})$ to space $L^{p}({\omega})$ with $1<p_1,p_2<\infty$ and $\frac{1}{p}=\frac{1}{p_1}+\frac{1}{p_2}$.
\end{lemma}

\section{Proof of the main theorems}
In this section, We will mainly give out the proof of Theorem 1.8. and Theorem 1.9.\\
{\bf Proof of Theorem 1.8.}  We just point out that the estimate for the strong type is almost the same as the weak type. Here we only present the proof of the strong type estimate.
For any fixed ball $B=B(x_0,r_B)\subset\mathcal{X}$ and $2B:=B(x_0,2r_B)$. We decompose $f_{j}$ as
$$f_j=f_{j}^{0}+f_{j}^{\infty}=f_{j}\chi_{2B}+f_{j}\chi_{\mathcal{X}\setminus2B}, j=1,2.$$

By the linearity of $T_{\theta}$ and the Minkowski inequality, we have
\begin{eqnarray*}
&&\frac{1}{\phi(\omega(B))^{\frac{1}{p}}}\bigg[\int_{B}|T_{\theta}(f_{1},f_{2})(x)|^{p}\omega(x)d\mu(x)\bigg]^{\frac{1}{p}}\\
&&\quad\quad\leq \frac{1}{\phi(\omega(B))^{\frac{1}{p}}}\bigg[\int_{B}|T_{\theta}(f^{0}_{1},f^{0}_{2})(x)|^{p}\omega(x)d\mu(x)\bigg]^{\frac{1}{p}}\\
&&\quad\quad\quad +\frac{1}{\phi(\omega(B))^{\frac{1}{p}}}\bigg[\int_{B}|T_{\theta}(f^{0}_{1},f^{\infty}_{2})(x)|^{p}\omega(x)d\mu(x)\bigg]^{\frac{1}{p}}\\
&&\quad\quad\quad +\frac{1}{\phi(\omega(B))^{\frac{1}{p}}}\bigg[\int_{B}|T_{\theta}(f^{\infty}_{1},f^{0}_{2})(x)|^{p}\omega(x)d\mu(x)\bigg]^{\frac{1}{p}}\\
&&\quad\quad\quad +\frac{1}{\phi(\omega(B))^{\frac{1}{p}}}\bigg[\int_{B}|T_{\theta}(f^{\infty}_{1},f^{\infty}_{2})(x)|^{p}\omega(x)d\mu(x)\bigg]^{\frac{1}{p}}\\
&&\quad\quad:=I_{1}+I_{2}+I_{3}+I_{4}.
\end{eqnarray*}
For the term $I_1$. By Lemma 2.2 and 2.5, we can deduce that
\begin{eqnarray*}
I_1&\leq&\frac{1}{\phi(\omega(B))^{\frac{1}{p}}}\bigg[\int_{\mathcal{X}}|T_{\theta}(f^{0}_{1},f^{0}_{2})(x)|^{p}\omega(x)d\mu(x)\bigg]^{\frac{1}{p}}\\
&\lesssim& \prod^{2}_{j=1}\|f_j\|_{\mathcal{M}^{p_j,\phi}(\omega)}\bigg[\frac{\phi({\omega}(2B)}{\phi({\omega}(B)}\bigg]^{\frac{1}{p}}\\
&\lesssim& \prod^{2}_{j=1}\|f_j\|_{\mathcal{M}^{p_j,\phi}(\omega)}\bigg[\frac{{\omega}(2B)}{{\omega}(B)}\bigg]^{\frac{1}{p}}\\
&\lesssim& \prod^{2}_{j=1}\|f_j\|_{\mathcal{M}^{p_j,\phi}(\omega)}.
\end{eqnarray*}
To estimates $I_2$, we first consider $|T_{\theta}(f^{0}_1,f^{\infty}_2)(x)|$, for any $x\in B $ and $y_1\in 2B, y_{2}\in (2B)^{c}$, then $d(x,y_{2})\sim d(x_0,y_2)$ and $V(x,y_2)\sim V(x_0,y_2)$.
By applying (1.4) and (1.6) with H\"{o}lder inequality, we can get
\begin{eqnarray*}
&&|T_{\theta}(f^{0}_1,f^{\infty}_2)(x)|\\
&&\quad\leq\int_{\mathcal{X}}\int_{\mathcal{X}}|K(x,y_1,y_2)||f^{0}_{1}(y_1)||f^{\infty}_2(y)|d\mu(y_1)d\mu(y_2)\\
&&\quad\lesssim\int_{2B}\int_{\mathcal{X}\backslash 2B}\frac{|f_{1}(y_1)||f_2(y)|}{{[\sum^{2}_{j=1}V(x,y_{j})]^{2}}}d\mu(y_1)d\mu(y_2)\\
&&\quad\lesssim\int_{2B}|f_{1}(y_1)|[\omega(y_1)]^{\frac{1}{p_1}-\frac{1}{p_1}}d\mu(y_1)\bigg\{\sum^{\infty}_{k=1}\int_{2^{k+1}B\backslash 2^{k}B}\frac{|f_2(y_2)|}{{[\sum^{2}_{j=1}V(x,y_{2})]^{2}}}d\mu(y_2)\bigg\}\\
&&\quad\lesssim [\phi(\omega(2B))]^{\frac{1}{p_1}}\bigg(\frac{1}{\phi(\omega(2B))}
\int_{2B}|f_{1}(y_1)|^{p_1}\omega(y_1)\bigg)^{\frac{1}{p_1}}\bigg(\int_{2B}[\omega(y_1)]^{-\frac{p^{\prime}_1}{p_1}}\bigg)^{\frac{1}{p^{\prime}_1}}\\
&&\quad\quad\times \bigg\{\sum^{\infty}_{k=1}\frac{1}{[\mu(2^{k}B)]^{2}}\int_{2^{k+1}B}|f_{2}(y_2)|[\omega(y_2)]^{\frac{1}{p_2}-\frac{1}{p_2}}d\mu(y_2)\bigg\}\\
&&\quad\lesssim \|f_1\|_{\mathcal{M}^{p_1,\phi}(\omega)}[\phi(\omega(2B))]^{\frac{1}{p_1}}
\bigg(\int_{2B}[\omega(y_1)]^{-\frac{p^{\prime}_1}{p_1}}\bigg)^{\frac{1}{p^{\prime}_1}}\bigg\{\sum^{\infty}_{k=1}\frac{1}{[\mu(2^{k}B)]^{2}}{\phi(\omega(2^{k+1}B))}^{\frac{1}{p_{2}}}\\
&&\quad\quad\times \bigg(\frac{1}{\phi(\omega(2^{k+1}B))}
\int_{2^{k+1}B}|f_{2}(y_2)|^{p_2}\omega(y_2)d\mu(y_2)\bigg)^{\frac{1}{p_2}}\bigg(\int_{2^{k+1}B}[\omega(y_{2})]^{-\frac{p^{\prime}_{2}}{p_2}}d\mu(y_{2})\bigg)^{\frac{1}{p^{\prime}_{2}}}\bigg\}\\
&&\quad\lesssim \prod^{2}_{j=1}\|f_j\|_{\mathcal{M}^{p_j,\phi}(\omega)}[\phi(\omega(2B))]^{\frac{1}{p_1}}\frac{\mu(2B)}{\omega(2B)^{\frac{1}{p_1}}}\bigg\{\sum^{\infty}_{k=1}
\frac{{\phi(\omega(2^{k+1}B))}^{\frac{1}{p_{2}}}}{[\mu(2^{k}B)]^{2}}\frac{\mu(2^{k+1}B)}{\omega(2^{k+1}B)^{\frac{1}{p_2}}}\bigg\}\\
&&\quad\lesssim \prod^{2}_{j=1}\|f_j\|_{\mathcal{M}^{p_j,\phi}(\omega_{j})}\bigg[\frac{\phi(\omega(B))}{\omega(B)}\bigg]^{\frac{1}{p}}.
\end{eqnarray*}

Furthermore, we can deduce that
\begin{eqnarray*}
I_2 &\leq& \frac{1}{\phi({\omega}(B))^{\frac{1}{p}}}\bigg[\int_{B}|T_{\theta}(f^{0}_{1},f^{\infty}_{2})(x)|^{p}{{\omega}}(x)d\mu(x)\bigg]^{\frac{1}{p}}\\
&\lesssim &\prod^{2}_{j=1}\|f_j\|_{\mathcal{M}^{p_j,\phi}(\omega)}.
\end{eqnarray*}
Since the estimates for and $I_2$ and $I_3$ are similar, hence
$$I_3\lesssim \prod^{2}_{j=1}\|f_j\|_{\mathcal{M}^{p_j,\phi}(\omega)}.$$

For the term $I_4$, we first estimate $T_{\theta}(f^{\infty}_{1},f^{\infty}_{2})(x)$ with $x\in B$. By (1.4) and  H\"{o}lder inequality, Definition 1.1 and 1.6, we obtain
\begin{eqnarray*}
\hspace{-10mm}
&&|T_{\theta}(f^{\infty}_1,f^{\infty}_2)(x)|\\
&&\quad\quad\leq\int_{\mathcal{X}\backslash 2B}\int_{\mathcal{X}\backslash2B}|K(x,y_1,y_2)||f_{1}(y_1)||f_2(y_{2})|d\mu(y_1)d\mu(y_2)\\
&&\quad\quad\lesssim\int_{\mathcal{X}\backslash 2B}\int_{\mathcal{X}\backslash2B}\frac{|f_{1}(y_1)||f_2(y)|}{{[\sum^{2}_{j=1}V(x,y_{j})]^{2}}}d\mu(y_1)d\mu(y_2)\\
&&\quad\quad\lesssim\prod^{2}_{j=1}\bigg\{\sum^{\infty}_{k=1}\int_{2^{k+1}B\backslash2^{k}B}\frac{|f_{j}(y_j)|}{V(x_0,y_j)}d\mu(y_j)\bigg\}\\
&&\quad\quad\lesssim\prod^{2}_{j=1}\bigg\{\sum^{\infty}_{k=1}\frac{1}{\mu(2^{k}B)}\int_{2^{k+1}B}|f_{j}(y_j)|[\omega(y_j)]^{\frac{1}{p_j}}[\omega(y_j)]^{-\frac{1}{p_j}}d\mu(y_j)\bigg\}\\
&&\quad\quad\lesssim\prod^{2}_{j=1}\bigg\{\sum^{\infty}_{k=1}\frac{{\phi(\omega(2^{k+1}B))}^{\frac{1}{p_j}}}{\mu(2^{k}B)}\bigg(\frac{1}{\phi(\omega(2^{k+1}B))}\int_{2^{k+1}B}|f_{j}(y_j)|^{p_{j}}\omega(y_j)d\mu{y_{j}}\bigg)^{\frac{1}{p_{j}}}\\
&&\quad\quad\quad\times \bigg(\int_{2^{k+1}B}[\omega(y_{j})]^{-\frac{p^{\prime}_j}{p_j}}d\mu(y_{j})\bigg)^{\frac{1}{p^{\prime}_j}}\\
&&\quad\quad\lesssim\prod^{2}_{j=1}\|f_j\|_{\mathcal{M}^{p_j,\phi}(\omega)}\bigg\{\sum^{\infty}_{k=1}\frac{{\phi(\omega(2^{k+1}B))}^{\frac{1}{p_j}}}{\mu(2^{k}B)}\cdot\frac{\mu(2^{k+1}B)}{(\omega(2^{k+1}B))^{\frac{1}{p_{j}}}}\bigg\}\\
&&\quad\quad\lesssim\prod^{2}_{j=1}\|f_j\|_{\mathcal{M}^{p_j,\phi}(\omega)}\bigg[\frac{\phi(\omega(B))}{\omega(B)}\bigg]^{\frac{1}{p}}.
\end{eqnarray*}
Furthermore, together with Definition 1.1, implies that
\begin{eqnarray*}
I_4 &\leq& \frac{1}{\phi({\omega}(B))^{\frac{1}{p}}}\bigg[\int_{B}|T_{\theta}(f^{\infty}_{1},f^{\infty}_{2})(x)|^{p}{{\omega}}(x)d\mu(x)\bigg]^{\frac{1}{p}}\\
&\lesssim &\prod^{2}_{j=1}\|f_j\|_{\mathcal{M}^{p_j,\phi}(\omega)}.
\end{eqnarray*}
Combining the estimate for $I_1, I_2, I_3$ and $I_4$, the theorem 1.8 is proved.\\
{\bf Proof of Theorem 1.9.}  We just point out that the estimate for the strong type is almost the same as the weak type. Here we omit the proof details of weak type estimate. For any fixed ball $B=B(x_0,r_B)\subset\mathcal{X}$ and $2B:=B(x_0,2r_B)$. We decompose $f_{j}$ as
$$f_j=f_{j}^{0}+f_{j}^{\infty}=f_{j}\chi_{2B}+f_{j}\chi_{\mathcal{X}\setminus2B}, j=1,2.$$

Then, write
\begin{eqnarray*}
&&\|[b_1,b_2,T_{\theta}](f_1,f_2)\|_{\mathcal{M}^{p,\phi}(\omega)}\\
&& \quad \leq\|[b_1,b_2,T_{\theta}](f^{0}_1,f^{0}_2)\|_{\mathcal{M}^{p,\phi}(\omega)}+\|[b_1,b_2,T_{\theta}](f^{0}_1,f^{\infty}_2)\|_{\mathcal{M}^{p,\phi}(\omega)}\\
&&\quad\quad +\|[b_1,b_2,T_{\theta}](f^{\infty}_1,f^{\infty}_2)\|_{\mathcal{M}^{p,\phi}(\omega)}+\|[b_1,b_2,T_{\theta}](f^{\infty}_1,f^{\infty}_2)\|_{\mathcal{M}^{p,\phi}(\omega)}\\
&&\quad:=J_1+J_2+J_3+J_4.
\end{eqnarray*}

The estimates for $J_{1}$ goes as follows. By applying Definition 1.1 and Lemma 2.5, we obtain that
\begin{eqnarray*}
J_{1}&=&\mathop{\sup}\limits_{B}\frac{1}{\phi(\omega(B))^{\frac{1}{p}}}\bigg\{\int_{B}|[b_1,b_2,T_{\theta}](f^{0}_{1},f^{0}_{2})(x)|^{p}\omega(x)d\mu(x)\bigg\}^{\frac{1}{p}}\\
&\leq& \mathop{\sup}\limits_{B}\frac{1}{\phi(\omega(B))^{\frac{1}{p}}}\|[b_1,b_2,T_{\theta}](f^{0}_{1},f^{0}_{2})\|_{L^{p}(\omega)}\\
&\lesssim&\|b_1\|_{BMO(\mu)}\|b_2\|_{BMO(\mu)}\mathop{\sup}\limits_{B}\frac{1}{\phi(\omega(B))^{\frac{1}{p}}}\bigg(\int_{2B}|f_{1}(y_{1})|^{p_1}\omega(y_1)d\mu(y_1)\bigg)^{\frac{1}{p_1}}\\
&&\times\bigg(\int_{2B}|f_{2}(y_{2})|^{p_2}\omega(y_2)d\mu(y_2)\bigg)^{\frac{1}{p_2}}\\
&\lesssim&\|b_1\|_{BMO(\mu)}\|b_2\|_{BMO(\mu)}\|f_{1}\|_{\mathcal{M}^{p_1,\phi}(\omega)}\|f_{2}\|_{\mathcal{M}^{p_2,\phi}(\omega)}\mathop{\sup}\limits_{B}\bigg(\frac{\phi(\omega(2B))}{\phi(\omega(B))}\bigg)^{\frac{1}{p}}\\
&\lesssim& \prod^{2}_{j=1}\|b_j\|_{BMO(\mu)}\|f_j\|_{\mathcal{M}^{p_j,\phi}(\omega)}\bigg(\frac{(\omega(2B))}{(\omega(B))}\bigg)^{\frac{1}{p}}\\
&\lesssim& \prod^{2}_{j=1}\|b_j\|_{BMO(\mu)}\|f_j\|_{\mathcal{M}^{p_j,\phi}(\omega)}.
\end{eqnarray*}

For any $x\in B$, write
\begin{eqnarray*}
&&|[b_1,b_2,T_{\theta}](f^{0}_{1},f^{\infty}_{2})(x)|\\
&&\quad\leq\int_{2B}\int_{\mathcal{X}\backslash 2B}|b_{1}(x)-b_1(y_1)||b_{2}(x)-b_2(y_2)||K(x,y_1,y_2)||f_1(y_1)||f_{2}(y_2)|d\mu(y_1)d\mu(y_2)\\
&&\quad\lesssim\int_{2B}\int_{\mathcal{X}\backslash 2B}\frac{|b_{1}(x)-b_1(y_1)||b_{2}(x)-b_2(y_2)|}{[\sum^{2}_{j=1}V(x,y_{j})]^{2}}|f_1(y_1)||f_{2}(y_2)|d\mu(y_1)d\mu(y_2)\\
&&\quad\lesssim\int_{2B}|b_{1}(x)-b_1(y_1)||f_{1}(y_{1})|d\mu(y_1)\bigg(\int_{\mathcal{X}\backslash 2B}\frac{|b_{2}(x)-b_2(y_2)|}{[V(x,y_{2})]^{2}}|f_{2}(y_2)|d\mu(y_2)\bigg)\\
&&\quad\lesssim |b_{1}(x)-(b_1)_{2B}|\int_{2B}|f_{1}(y_{1})|d\mu(y_1)\bigg(\int_{\mathcal{X}\backslash 2B}\frac{|b_{2}(x)-b_2(y_2)|}{[V(x,y_{2})]^{2}}|f_{2}(y_2)|d\mu(y_2)\bigg)\\
&&\quad\quad +\int_{2B}|b_{1}(y_1)-(b_1)_{2B}||f_{1}(y_{1})|d\mu(y_1)\bigg(\int_{\mathcal{X}\backslash 2B}\frac{|b_{2}(x)-b_2(y_2)|}{[V(x,y_{2})]^{2}}|f_{2}(y_2)|d\mu(y_2)\bigg)\\
&&\quad:=J_{21}+J_{22}.
\end{eqnarray*}
By virtue of H\"{o}lder's inequality and lemma 2.2-2.4, it follows that
\begin{eqnarray*}
J_{21}&=&|b_{1}(x)-(b_1)_{2B}|\int_{2B}|f_{1}(y_{1})|d\mu(y_1)\bigg(\int_{\mathcal{X}\backslash 2B}\frac{|b_{2}(x)-b_2(y_2)|}{[V(x,y_{2})]^{2}}|f_{2}(y_2)|d\mu(y_2)\bigg)\\
&\leq&|b_{1}(x)-(b_1)_{2B}|\bigg(\int_{2B}|f_{1}(y_1)|^{p_1}\omega(y_1)d\mu(y_1)\bigg)^{\frac{1}{p_1}}\bigg(\int_{2B}\omega(y_1)^{-\frac{p^{\prime}_1}{p_1}}d\mu(y_1)\bigg)^{\frac{1}{p^{\prime}_1}}\\
&&\times\bigg(\sum^{\infty}_{k=1}\int_{2^{k+1}B\backslash2^{k}B}\frac{|b_{2}(x)-(b_2)_{2B}+(b_2)_{2B}-b_2(y_2)|}{[V(x,y_{2})]^{2}}|f_{2}(y_2)|d\mu(y_2)\bigg)\\
&\leq&|b_{1}(x)-(b_1)_{2B}|\phi(\omega(2B))^{\frac{1}{p_1}}\|f_1\|_{\mathcal{M}^{p_1,\phi}(\omega)}\frac{\mu(2B)}{\omega(2B)^{\frac{1}{p_1}}}\\
&&\times\bigg(|b_{2}(x)-(b_2)_{2B}|\sum^{\infty}_{k=1}\int_{2^{k+1}B\backslash2^{k}B}\frac{|f_{2}(y_2)|}{[V(x,y_{2})]^{2}}|d\mu(y_2)\bigg)\\
&&+|b_{1}(x)-(b_1)_{2B}|\phi(\omega(2B))^{\frac{1}{p_1}}\|f_1\|_{\mathcal{M}^{p_1,\phi}(\omega)}\frac{\mu(2B)}{\omega(2B)^{\frac{1}{p_1}}}\\
&&\times\bigg(\sum^{\infty}_{k=1}\int_{2^{k+1}B\backslash2^{k}B}\frac{|(b_2)_{2B}-b_2(y_2)|}{[V(x,y_{2})]^{2}}||f_2(y_2)|d\mu(y_2)\bigg)\\
&\leq&|b_{1}(x)-(b_1)_{2B}|\phi(\omega(2B))^{\frac{1}{p_1}}\|f_1\|_{\mathcal{M}^{p_1,\phi}(\omega)}\frac{\mu(2B)}{\omega(2B)^{\frac{1}{p_1}}}\\
&&\times\bigg\{|b_{2}(x)-(b_2)_{2B}|\sum^{\infty}_{k=1}\frac{1}{[\mu(2^{k}B)]^{2}}\bigg(\int_{2^{k+1}B}|f_{2}(y_2)|^{p_2}\omega(y_2)d\mu{(y_2)}
\bigg)^{\frac{1}{p_2}}\cdot\frac{\mu(2^{k+1}B)}{\omega(2^{k+1}B)^{\frac{1}{p_{2}}}}\bigg\}\\
&&+|b_{1}(x)-(b_1)_{2B}|\phi(\omega(2B))^{\frac{1}{p_1}}\|f_1\|_{\mathcal{M}^{p_1,\phi}(\omega)}\frac{\mu(2B)}{\omega(2B)^{\frac{1}{p_1}}}\\
&&\times\bigg\{\sum^{\infty}_{k=1}\frac{|(b_2)_{2B}-(b_2)_{2^{k+1}B}|}{[\mu(2^{k}B)]^{2}}\int_{2^{k+1}B}|f_{2}(y_{2})|d\mu(y_2)\\
&&+\sum^{\infty}_{k=1}\frac{1}{[\mu(2^{k}B)]^{2}}\int_{2^{k+1}B}|b_{2}(y_2)-(b_2)_{2^{k+1}B}||f_{2}(y_{2})|d\mu(y_2)\bigg\}\\
&\leq& |b_{1}(x)-(b_1)_{2B}||b_{2}(x)-(b_2)_{2B}|\bigg[\frac{\phi(\omega(B))}{\omega(B)}\bigg]^{\frac{1}{p}}\|f_1\|_{\mathcal{M}^{p_1,\phi}(\omega)}\|f_2\|_{\mathcal{M}^{p_2,\phi}(\omega)}\\
&&+|b_{1}(x)-(b_1)_{2B}|\bigg[\frac{\phi(\omega(B))}{\omega(B)}\bigg]^{\frac{1}{p_1}}\|f_1\|_{\mathcal{M}^{p_1,\phi}(\omega)}\\
&&\times\bigg\{\sum^{\infty}_{k=1}(k+1)\|b_{2}\|_{BMO(\mu)}\bigg[\frac{\phi(\omega(2^{k+1}B))}{\omega(2^{k+1}B)}\bigg]^{\frac{1}{p_2}}\|f_2\|_{\mathcal{M}^{p_2,\phi}(\omega)}\bigg\}\\
&&+|b_{1}(x)-(b_1)_{2B}|\phi(\omega(2B))^{\frac{1}{p_1}}\|f_1\|_{\mathcal{M}^{p_1,\phi}(\omega)}\frac{\mu(2B)}{\omega(2B)^{\frac{1}{p_1}}}\\
&&\times\sum^{\infty}_{k=1}\frac{1}{[\mu(2^{k}B)]^{2}}\bigg(\int_{2^{k+1}B}|f_{2}(y_{2})|^{p_2}\omega(y_2)d\mu(y_2)\bigg)^{\frac{1}{p_2}}\\
&&\times\bigg(\int_{2^{k+1}B}|b_{2}(y_2)-(b_2)_{2^{k+1}B}|^{p^{\prime}_2}[\omega(y_2)]^{-\frac{p^{\prime}_2}{p_2}}d\mu(y_2)\bigg)^{\frac{1}{p^{\prime}_2}}\\
&\leq& |b_{1}(x)-(b_1)_{2B}||b_{2}(x)-(b_2)_{2B}|\bigg[\frac{\phi(\omega(B))}{\omega(B)}\bigg]^{\frac{1}{p}}\|f_1\|_{\mathcal{M}^{p_1,\phi}(\omega)}\|f_2\|_{\mathcal{M}^{p_2,\phi}(\omega)}\\
&&+|b_{1}(x)-(b_1)_{2B}|\bigg[\frac{\phi(\omega(B))}{\omega(B)}\bigg]^{\frac{1}{p_1}}\|f_1\|_{\mathcal{M}^{p_1,\phi}(\omega)}\\
&&\times\bigg\{\sum^{\infty}_{k=1}(k+1)\|b_{2}\|_{BMO(\mu)}\bigg[\frac{\phi(\omega(2^{k+1}B))}{\omega(2^{k+1}B)}\bigg]^{\frac{1}{p_2}}\|f_2\|_{\mathcal{M}^{p_2,\phi}(\omega)}\bigg\}\\
&&+|b_{1}(x)-(b_1)_{2B}|\phi(\omega(2B))^{\frac{1}{p_1}}\|f_1\|_{\mathcal{M}^{p_1,\phi}(\omega)}\frac{\mu(2B)}{\omega(2B)^{\frac{1}{p_1}}}\\
&&\times\sum^{\infty}_{k=1}\frac{1}{[\mu(2^{k}B)]^{2}}\bigg(\int_{2^{k+1}B}|f_{2}(y_{2})|^{p_2}\omega(y_2)d\mu(y_2)\bigg)^{\frac{1}{p_2}}\|b_2\|_{BMO(\mu)}\frac{\mu(2^{k+1}B)}{[\omega(2^{k+1}B)]^{\frac{1}{p_2}}}\\
&\lesssim & |b_{1}(x)-(b_1)_{2B}||b_{2}(x)-(b_2)_{2B}|\|f_1\|_{\mathcal{M}^{p_1,\phi}(\omega)}\|f_2\|_{\mathcal{M}^{p_2,\phi}(\omega)}\bigg[\frac{\phi(\omega(B))}{\omega(B)}\bigg]^{\frac{1}{p}}\\
&&+|b_{1}(x)-(b_1)_{2B}|\|b_{2}\|_{BMO(\mu)}\|f_1\|_{\mathcal{M}^{p_1,\phi}(\omega)}\|f_2\|_{\mathcal{M}^{p_2,\phi}(\omega)}\bigg[\frac{\phi(\omega(B))}{\omega(B)}\bigg]^{\frac{1}{p}}.\\
\end{eqnarray*}

As for the term $J_{22}$, by applying H\"{o}lder's inequality and lemma 2.3, we have
\begin{eqnarray*}
J_{22}&=&\int_{2B}|b_{1}(y_1)-(b_1)_{2B}||f_{1}(y_{1})|d\mu(y_1)\bigg(\int_{\mathcal{X}\backslash 2B}\frac{|b_{2}(x)-b_2(y_2)|}{[V(x,y_{2})]^{2}}|f_{2}(y_2)|d\mu(y_2)\bigg)\\
&\leq&\bigg(\int_{2B}|f_{1}(y_1)|^{p_1}\omega(y_1)d\mu(y_1)\bigg)^{\frac{1}{p_1}}\bigg(\int_{2B}|b_1(y_1)-(b_1)_{2B}|^{p^{\prime}_1}[\omega(y_{1})]^{-\frac{p^{\prime}_1}{p_1}}d\mu(y_1)\bigg)^{\frac{1}{p^{\prime}_1}}\\
&&\times\bigg\{|b_{2}(x)-(b_2)_{2B}|\int_{\mathcal{X}\backslash 2B}\frac{|f_{2}(y_2)|}{[V(x,y_{2})]^{2}}d\mu(y_2)\\
&&+\int_{\mathcal{X}\backslash 2B}\frac{|b_2(y_2)-(b_2)_{2B}|}{[V(x,y_{2})]^{2}}|f_{2}(y_2)|d\mu(y_2)\bigg\}\\
&\lesssim & \|b_{1}\|_{BMO(\mu)}\|f_{1}\|_{\mathcal{M}^{p_1,\phi}(\omega)}[\phi(\omega(2B))]^{\frac{1}{p_1}}\frac{\mu(2B)}{[\omega(2B)]^{\frac{1}{p_1}}}\\
&&\times\bigg\{|b_{2}(x)-(b_2)_{2B}|\sum^{\infty}_{k=1}\int_{2^{k+1}B\backslash2^{k}B}\frac{|f_{2}(y_2)|}{[V(x,y_{2})]^{2}}d\mu(y_2)\\
&&+\sum^{\infty}_{k=1}\int_{2^{k+1}B\backslash2^{k}B}\frac{|b_2(y_2)-(b_2)_{2B}|}{[V(x,y_{2})]^{2}}|f_{2}(y_2)|d\mu(y_2)\bigg\}\\
&\lesssim & \|b_{1}\|_{BMO(\mu)}\|f_{1}\|_{\mathcal{M}^{p_1,\phi}(\omega)}[\phi(\omega(2B))]^{\frac{1}{p_1}}\frac{\mu(2B)}{[\omega(2B)]^{\frac{1}{p_1}}}\\
&&\times\bigg\{|b_{2}(x)-(b_2)_{2B}|\sum^{\infty}_{k=1}\frac{1}{[\mu(2^{k}B)]^{2}}\int_{2^{k+1}B}|f_{2}(y_2)|d\mu(y_2)\\
&&+\sum^{\infty}_{k=1}\frac{1}{[\mu(2^{k}B)]^{2}}\int_{2^{k+1}B}|b_2(y_2)-(b_2)_{2B}||f_{2}(y_2)|d\mu(y_2)\bigg\}\\
&\lesssim & \|b_{1}\|_{BMO(\mu)}\|f_{1}\|_{\mathcal{M}^{p_1,\phi}(\omega)}[\phi(\omega(2B))]^{\frac{1}{p_1}}\frac{\mu(2B)}{[\omega(2B)]^{\frac{1}{p_1}}}\\
&&\times\bigg\{|b_{2}(x)-(b_2)_{2B}|\sum^{\infty}_{k=1}\frac{1}{[\mu(2^{k}B)]^{2}}\bigg(\int_{2^{k+1}B}|f_{2}(y_2)|^{p_2}\omega(y_2)d\mu(y_2)\bigg)^{\frac{1}{p_2}}\frac{\mu(2^{k+1}B)}{[\omega(2^{k+1}B)]^{\frac{1}{p_2}}}\\
&&+\sum^{\infty}_{k=1}\frac{|(b_2)_{2B}-(b_2)_{2^{k+}B}|}{[\mu(2^{k}B)]^{2}}\bigg(\int_{2^{k+1}B}|f_{2}(y_2)|^{p_2}\omega(y_2)d\mu(y_2)\bigg)^{\frac{1}{p_2}}\frac{\mu(2^{k+1}B)}{[\omega(2^{k+1}B)]^{\frac{1}{p_2}}}\\
&&+\sum^{\infty}_{k=1}\frac{1}{[\mu(2^{k}B)]^{2}}\bigg(\int_{2^{k+1}B}|f_{2}(y_2)|^{p_2}\omega(y_2)d\mu(y_2)\bigg)^{\frac{1}{p_2}}\\
&&\times\bigg(\int_{2^{k+1}B}|b_{2}(x)-(b_{2})_{2^{k+1}B}|^{p^{\prime}_2}[\omega(y_2)]^{-\frac{p^{\prime}_2}{p_2}}d\mu(y_2)\bigg)^{\frac{1}{p^{\prime}_2}}\bigg\}\\
&\lesssim & \|b_{1}\|_{BMO(\mu)}\|f_{1}\|_{\mathcal{M}^{p_1,\phi}(\omega)}[\phi(\omega(2B))]^{\frac{1}{p_1}}\frac{\mu(2B)}{[\omega(2B)]^{\frac{1}{p_1}}}\\
&&\times\bigg\{\|f_2\|_{\mathcal{M}^{p_2,\phi}(\omega)}|b_{2}(x)-(b_2)_{2B}|\sum^{\infty}_{k=1}\bigg[\frac{\phi(\omega(2^{k+1}B))}{\omega(2^{k+1}B)}\bigg]^{\frac{1}{p_2}}\\
&&+\|b_{2}\|_{BMO(\mu)}\|f_2\|_{\mathcal{M}^{p_2,\phi}(\omega)}\sum^{\infty}_{k=1}(k+1)\bigg[\frac{\phi(\omega(2^{k+1}B))}{\omega(2^{k+1}B)}\bigg]^{\frac{1}{p_2}}\\
&\lesssim & (|b_{2}(x)-(b_2)_{2B}|+\|b_2\|_{BMO(\mu)})\|b_{1}\|_{BMO(\mu)}\|f_{1}\|_{\mathcal{M}^{p_1,\phi}(\omega)}\|f_2\|_{\mathcal{M}^{p_2,\phi}(\omega)}\bigg[\frac{\phi(\omega(B))}{\omega(B)}\bigg]^{\frac{1}{p}},
\end{eqnarray*}

which, together with the estimates of $J_{21}$ and lemma 2.4, we can deduce that
\begin{eqnarray*}
&&\|[b_1,b_2,T_{\theta}](f^{0}_1,f^{\infty}_2)\|_{\mathcal{M}^{p,\phi}(\omega)}\\
&&\quad=\mathop{\sup}\limits_{B}\frac{1}{\phi(\omega(B))^{\frac{1}{p}}}\bigg(\int_{B}|[b_1,b_2,T_{\theta}](f^{0}_1,f^{\infty}_2)(x)|^{p}\omega(x)d\mu(x)\bigg)^{\frac{1}{p}}\\
&&\quad\lesssim \|b_1\|_{BMO(\mu)}\|b_2\|_{BMO(\mu)}\|f_{1}\|_{\mathcal{M}^{p_1,\phi}(\omega)}\|f_2\|_{\mathcal{M}^{p_2,\phi}(\omega)}
\bigg[\frac{\omega(B)}{\phi(\omega(B))}\bigg]^{\frac{1}{p}}\bigg[\frac{\phi(\omega(B))}{\omega(B)}\bigg]^{\frac{1}{p}}\\
&&\quad\quad+\|b_{2}\|_{BMO(\mu)}\|f_{1}\|_{\mathcal{M}^{p_1,\phi}(\omega)}\|f_2\|_{\mathcal{M}^{p_2,\phi}(\omega)}\bigg[\frac{\phi(\omega(B))}{\omega(B)}\bigg]^{\frac{1}{p}}\\
&&\quad\quad\times\mathop{\sup}\limits_{B}\frac{1}{\phi(\omega(B))^{\frac{1}{p}}}\bigg(\int_{B}|b_1(x)-(b_1)_{2B}|^{p}\omega{(x)}d\mu(x)\bigg)^{\frac{1}{p}}\\
&&\quad\quad+\|b_{1}\|_{BMO(\mu)}\|f_{1}\|_{\mathcal{M}^{p_1,\phi}(\omega)}\|f_2\|_{\mathcal{M}^{p_2,\phi}(\omega)}\bigg[\frac{\phi(\omega(B))}{\omega(B)}\bigg]^{\frac{1}{p}}\\
&&\quad\quad\times\mathop{\sup}\limits_{B}\frac{1}{\phi(\omega(B))^{\frac{1}{p}}}\bigg(\int_{B}|b_2(x)-(b_2)_{2B}|^{p}\omega{(x)}d\mu(x)\bigg)^{\frac{1}{p}}\\
&&\quad\quad+\|f_{1}\|_{\mathcal{M}^{p_1,\phi}(\omega)}\|f_2\|_{\mathcal{M}^{p_2,\phi}(\omega)}\bigg[\frac{\phi(\omega(B))}{\omega(B)}\bigg]^{\frac{1}{p}}\\
&&\quad\quad\times\mathop{\sup}\limits_{B}\frac{1}{\phi(\omega(B))^{\frac{1}{p}}}\bigg(\int_{B}|b_1(x)-(b_1)_{2B}|^{p}|b_2(x)-(b_2)_{2B}|^{p}\omega{(x)}d\mu(x)\bigg)^{\frac{1}{p}}\\
&&\quad\lesssim \|b_1\|_{BMO(\mu)}\|b_2\|_{BMO(\mu)}\|f_{1}\|_{\mathcal{M}^{p_1,\phi}(\omega)}\|f_2\|_{\mathcal{M}^{p_2,\phi}(\omega)}\\
&&\quad\quad+\|f_{1}\|_{\mathcal{M}^{p_1,\phi}(\omega)}\|f_2\|_{\mathcal{M}^{p_2,\phi}(\omega)}\mathop{\sup}\limits_{B}\bigg[\frac{\phi(\omega(B))}{\omega(B)}\bigg]^{\frac{1}{p}}\frac{1}{\phi(\omega(B))^{\frac{1}{p}}}\\
&&\quad\quad\times\bigg(\int_{B}|b_1(x)-(b_1)_{2B}|^{p_{1}}\omega(x)d\mu(x)\bigg)^{\frac{1}{p_1}}\bigg(\int_{B}|b_2(x)-(b_2)_{2B}|^{p_{2}}\omega{(x)}d\mu(x)\bigg)^{\frac{1}{p_2}}\\
&&\quad\lesssim \|b_1\|_{BMO(\mu)}\|b_2\|_{BMO(\mu)}\|f_{1}\|_{\mathcal{M}^{p_1,\phi}(\omega)}\|f_2\|_{\mathcal{M}^{p_2,\phi}(\omega)}.
\end{eqnarray*}

Similarly, we get
$$J_3\lesssim \|b_1\|_{BMO(\mu)}\|b_2\|_{BMO(\mu)}\|f_{1}\|_{\mathcal{M}^{p_1,\phi}(\omega)}\|f_2\|_{\mathcal{M}^{p_2,\phi}(\omega)}.$$

For any $x\in B$, we write
\begin{eqnarray*}
&&|[b_1,b_2,T_{\theta}](f^{\infty}_{1},f^{\infty}_{2})(x)|\\
&&\quad\lesssim\int_{\mathcal{X}\backslash 2B}\int_{\mathcal{X}\backslash 2B}\frac{|b_{1}(x)-b_1(y_1)||b_{2}(x)-b_2(y_2)|}{[\sum^{2}_{j=1}V(x,y_{j})]^{2}}|f_1(y_1)||f_{2}(y_2)|d\mu(y_1)d\mu(y_2)\\
&&\quad\lesssim |b_1(x)-(b_1)_{2B}||b_2(x)-(b_2)_{2B}|\int_{\mathcal{X}\backslash 2B}\int_{\mathcal{X}\backslash 2B}\frac{|f_1(y_1)||f_{2}(y_2)|}{[\sum^{2}_{j=1}V(x,y_{j})]^{2}}d\mu(y_1)d\mu(y_2)\\
&&\quad\quad +|b_1(x)-(b_1)_{2B}|\int_{\mathcal{X}\backslash 2B}\int_{\mathcal{X}\backslash 2B}\frac{|b_2(y_{2})-(b_2)_{2B}|}{[\sum^{2}_{j=1}V(x,y_{j})]^{2}}|f_1(y_1)||f_{2}(y_2)|d\mu(y_1)d\mu(y_2)\\
&&\quad\quad +|b_2(x)-(b_2)_{2B}|\int_{\mathcal{X}\backslash 2B}\int_{\mathcal{X}\backslash 2B}\frac{|b_1(y_{1})-(b_1)_{2B}|}{[\sum^{2}_{j=1}V(x,y_{j})]^{2}}|f_1(y_1)||f_{2}(y_2)|d\mu(y_1)d\mu(y_2)\\
&&\quad\quad +\int_{\mathcal{X}\backslash 2B}\int_{\mathcal{X}\backslash 2B}\frac{|(b_1)_{2B}-b_1(y_1)||(b_{2})_{2B}-b_2(y_2)|}{[\sum^{2}_{j=1}V(x,y_{j})]^{2}}|f_1(y_1)||f_{2}(y_2)|d\mu(y_1)d\mu(y_2)\\
&&\quad:=J_{41}+J_{42}+J_{43}+J_{44}.
\end{eqnarray*}

By applying H\"{o}lder's inequality and Lemma 2.3, we deduce that
\begin{eqnarray*}
J_{41}&=&\prod^{2}_{i=1}|b_{i}(x)-(b_{i})_{2B}|\bigg\{\sum^{\infty}_{k=1}\int_{2^{k+1}B\backslash2^{k}B}|f_{1}(y_1)|\\
&&\times\bigg(\sum^{\infty}_{j=1}\int_{2^{j+1}B\backslash2^{j}B}\frac{|f_{2}(y_2)|}{[V(x,y_{1})+V(x,y_{2})]^{2}}d\mu(y_2)\bigg)d\mu(y_1)\bigg\}\\
&=&\prod^{2}_{i=1}|b_{i}(x)-(b_{i})_{2B}|\bigg\{\sum^{\infty}_{k=1}\int_{2^{k+1}B\backslash2^{k}B}|f_{1}(y_1)|\\
&&\times\bigg(\sum^{k}_{j=1}\int_{2^{j+1}B\backslash2^{j}B}\frac{|f_{2}(y_2)|}{[V(x,y_{1})+V(x,y_{2})]^{2}}d\mu(y_2)\\
&&+\sum^{\infty}_{j=k+1}\int_{2^{j+1}B\backslash2^{j}B}\frac{|f_{2}(y_2)|}{[V(x,y_{1})+V(x,y_{2})]^{2}}d\mu(y_2)\bigg)d\mu(y_1)\bigg\}\\
&\lesssim&\prod^{2}_{i=1}|b_{i}(x)-(b_{i})_{2B}|\bigg\{\sum^{\infty}_{k=1}\int_{2^{k+1}B\backslash2^{k}B}\frac{|f_{1}(y_1)|}{[V(x,y_{1})]^{2}}\bigg(\sum^{k}_{j=1}\int_{2^{j+1}B\backslash2^{j}B}|f_{2}(y_2)|d\mu(y_2)\bigg)\\
&&+\sum^{\infty}_{k=1}\int_{2^{k+1}B\backslash2^{k}B}|f_1(y_1)|\bigg(\sum^{\infty}_{j=k+1}\int_{2^{j+1}B\backslash2^{j}B}\frac{|f_2(y_2)|}{[V(x,y_{2})]^{2}}d\mu(y_2)\bigg)d\mu(y_1)\bigg\}\\
&\lesssim&\prod^{2}_{i=1}|b_{i}(x)-(b_{i})_{2B}|\bigg\{\sum^{\infty}_{k=1}\int_{2^{k+1}B\backslash2^{k}B}\frac{|f_{1}(y_1)|}{[V(x,y_{1})]^{2}}\bigg(\int_{2^{k+1}B}|f_{2}(y_2)|d\mu(y_2)\bigg)d\mu(y_1)\\
&&+\sum^{\infty}_{j=1}\int_{2^{j+1}B\backslash2^{j}B}\frac{|f_{2}(y_2)|}{[V(x,y_{2})]^{2}}\bigg(\sum^{j-1}_{k=1}\int_{2^{k+1}B\backslash2^{k}B}|f_{1}(y_1)|d\mu(y_1)\bigg)d\mu(y_2)\bigg\}\\
&\lesssim&\prod^{2}_{i=1}|b_{i}(x)-(b_{i})_{2B}|\bigg\{\sum^{\infty}_{k=1}\frac{1}{[\mu(2^{k}B)]^{2}}\int_{2^{k+1}B}|f_{1}(y_1)|\bigg(\int_{2^{k+1}B}|f_{2}(y_2)|d\mu(y_2)\bigg)d\mu(y_1)\\
&&+\sum^{\infty}_{j=1}\frac{1}{[\mu(2^{j}B)]^{2}}\int_{2^{j+1}B}|f_{2}(y_2)|\bigg(\int_{2^{j}B}|f_{1}(y_1)|d\mu(y_1)\bigg)d\mu(y_2)\bigg\}\\
&\lesssim&\prod^{2}_{i=1}|b_{i}(x)-(b_{i})_{2B}|\bigg\{\sum^{\infty}_{k=1}\frac{1}{[\mu(2^{k}B)]^{2}}\int_{2^{k+1}B}|f_{i}(y_i)|d\mu(y_i)\bigg\}\\
&\lesssim&\prod^{2}_{i=1}|b_{i}(x)-(b_{i})_{2B}|\bigg\{\sum^{\infty}_{k=1}\bigg(\int_{2^{k+1}B}|f_i(y_{i})|^{p_i}\omega(y_i)d\mu(y_i)\bigg)^{\frac{1}{p_i}}\frac{1}{[\omega(2^{k+1}B)]^{\frac{1}{p_i}}}\bigg\}\\
&\lesssim&\prod^{2}_{i=1}|b_{i}(x)-(b_{i})_{2B}|\|f_{i}\|_{\mathcal{M}^{p_i,\phi}(\omega)}\bigg\{\sum^{\infty}_{k=1}\bigg[\frac{\phi(\omega(2^{k+1}B))}{\omega(2^{k+1}B)}\bigg]^{\frac{1}{p_i}}\bigg\}\\
&\lesssim&\prod^{2}_{i=1}|b_{i}(x)-(b_{i})_{2B}|\|f_{i}\|_{\mathcal{M}^{p_i,\phi}(\omega)}\bigg[\frac{\phi(\omega(B))}{\omega(B)}\bigg]^{\frac{1}{p_i}}.\\
\end{eqnarray*}
For any $x\in B$, by Definition 1.1 , H\"{o}lder's inequality and Lemma 2.3, we can get
\begin{eqnarray*}
J_{42}&=&|b_1(x)-(b_1)_{2B}|\bigg\{\sum^{\infty}_{k=1}\int_{2^{k+1}B\backslash2^{k}B}|f_1(y_1)|\\
&&\times\bigg(\sum^{\infty}_{j=1}\int_{2^{j+1}B\backslash2^{j}B}\frac{|b_2(y_2)-(b_2)_{2B}|}{[V(x,y_{1})+V(x,y_{2})]^{2}}|f_{2}(y_2)|d\mu(y_2)\bigg)d\mu(y_1)\\
&\leq&|b_1(x)-(b_1)_{2B}|\bigg\{\sum^{\infty}_{k=1}\int_{2^{k+1}B\backslash2^{k}B}\frac{|f_1(y_1)|}{[V(x,y_{1})]^{2}}d\mu(y_1)\\
&&\times\bigg(\sum^{k}_{j=1}\int_{2^{j+1}B\backslash2^{j}B}|b_2(y_2)-(b_2)_{2B}||f_{2}(y_2)|d\mu(y_2)\bigg)\bigg\}\\
&&+|b_1(x)-(b_1)_{2B}|\bigg\{\sum^{\infty}_{k=1}\int_{2^{k+1}B\backslash2^{k}B}|f_1(y_1)|d\mu(y_1)\\
&&\times\bigg(\sum^{\infty}_{j=k+1}\int_{2^{j+1}B\backslash2^{j}B}\frac{|b_2(y_2)-(b_2)_{2B}|}{[V(x,y_{2})]^{2}}|f_{2}(y_2)|d\mu(y_2)\bigg)\bigg\}\\
&\leq&|b_1(x)-(b_1)_{2B}|\bigg\{\sum^{\infty}_{k=1}\int_{2^{k+1}B\backslash2^{k}B}\frac{|f_1(y_1)|}{[V(x,y_{1})]^{2}}d\mu(y_1)\\
&&\times\bigg(\int_{2^{k+1}B}|b_2(y_2)-(b_2)_{2B}||f_{2}(y_2)|d\mu(y_2)\bigg)\bigg\}\\
&&+|b_1(x)-(b_1)_{2B}|\bigg\{\sum^{\infty}_{j=1}\int_{2^{j+1}B\backslash2^{j}B}\frac{|b_2(y_2)-(b_2)_{2B}|}{[V(x,y_{2})]^{2}}|f_{2}(y_2)|d\mu(y_2)\\
&&\times\bigg(\sum^{j-1}_{k=1}\int_{2^{k+1}B\backslash2^{k}B}|f_1(y_1)|d\mu(y_1)\bigg)\bigg\}\\
&\leq&|b_1(x)-(b_1)_{2B}|\bigg\{\sum^{\infty}_{k=1}\frac{1}{[\mu(2^{k}B)]^{2}}\bigg(\int_{2^{k+1}B}|f_1(y_1)|^{p_1}\omega(y_1)d\mu(y_1)\bigg)^{\frac{1}{p_1}}\frac{\mu(2^{k+1}B)}{[\omega(2^{k+1}B)]^{\frac{1}{p_1}}}\\
&&\times\bigg(\int_{2^{k+1}B}|f_2(y_2)|^{p_2}\omega(y_2)d\mu(y_2)\bigg)^{\frac{1}{p_2}}\bigg(\int_{2^{k+1}B}|b_2(y_2)-(b_2)_{2B}|^{p^{\prime}_2}[\omega(y_2)]^{-\frac{p^{\prime}_2}{p_2}}d\mu(y_2)\bigg)^{\frac{1}{p^{\prime}_2}}\bigg\}\\
&&+|b_1(x)-(b_1)_{2B}|\bigg\{\sum^{\infty}_{j=1}\frac{1}{[\mu(2^{j}B)]^{2}}\bigg(\int_{2^{j+1}B}|f_2(y_2)|^{p_2}\omega(y_2)d\mu(y_2)\bigg)^{\frac{1}{p_2}}\\
&&\times\bigg(\int_{2^{j+1}B}|b_2(y_2)-(b_2)_{2B}|^{p^{\prime}_2}[\omega(y_2)]^{-\frac{p^{\prime}_2}{p_2}}d\mu(y_2)\bigg)^{\frac{1}{p^{\prime}_2}}\\
&&\times\bigg(\int_{2^{j}B}|f_1(y_1)|^{p_1}\omega(y_1)d\mu(y_1)\bigg)^{\frac{1}{p_1}}\frac{\mu(2^{j}B)}{[\omega(2^{j}B)]^{\frac{1}{p_1}}}\bigg\}\\
&\lesssim& |b_1(x)-(b_1)_{2B}|\|b_2\|_{BMO(\mu)}\|f_1\|_{\mathcal{M}^{p_1,\phi}(\omega)}\|f_2\|_{\mathcal{M}^{p_2,\phi}(\omega)}\bigg[\frac{\phi(\omega(B))}{\omega(B)}\bigg]^{\frac{1}{p}}.
\end{eqnarray*}

Similarty, we also get
$$J_{43}\lesssim |b_2(x)-(b_2)_{2B}|\|b_1\|_{BMO(\mu)}\|f_1\|_{\mathcal{M}^{p_1,\phi}(\omega)}\|f_2\|_{\mathcal{M}^{p_2,\phi}(\omega)}\bigg[\frac{\phi(\omega(B))}{\omega(B)}\bigg]^{\frac{1}{p}}.$$

For the last term $J_{44}$, by H\"{o}lder's inequality and Lemma 2.3, it follows that
\begin{eqnarray*}
J_{44}&=&\sum^{\infty}_{k=1}\int_{2^{k+1}B\backslash2^{k}B}|(b_1)_{2B}-b_1(y_1)||f_1(y_1)|\\
&&\times \bigg(\sum^{\infty}_{j=1}\int_{2^{j+1}B\backslash2^{j}B}\frac{|(b_{2})_{2B}-b_2(y_2)|}{[V(x,y_{1})+V(x,y_{2})]^{2}}|f_{2}(y_2)|d\mu(y_2)\bigg)d\mu(y_1)\\
&\leq&\sum^{\infty}_{k=1}\int_{2^{k+1}B\backslash2^{k}B}\frac{|(b_1)_{2B}-b_1(y_1)||f_1(y_1)|}{[V(x,y_{1})]^{2}}\\
&&\times \bigg(\sum^{k}_{j=1}\int_{2^{j+1}B\backslash2^{j}B}|(b_{2})_{2B}-b_2(y_2)||f_{2}(y_2)|d\mu(y_2)\bigg)d\mu(y_1)\\
&&+\sum^{\infty}_{k=1}\int_{2^{k+1}B\backslash2^{k}B}|(b_1)_{2B}-b_1(y_1)||f_1(y_1)|\\
&&\times \bigg(\sum^{\infty}_{j=k+1}\int_{2^{j+1}B\backslash2^{j}B}\frac{|(b_{2})_{2B}-b_2(y_2)||f_{2}(y_2)|}{[V(x,y_{2})]^{2}}d\mu(y_2)\bigg)d\mu(y_1)\\
&\leq& \sum^{\infty}_{k=1}\int_{2^{k+1}B\backslash2^{k}B}\frac{|(b_1)_{2B}-b_1(y_1)||f_1(y_1)|}{[V(x,y_{1})]^{2}}\\
&&\times \bigg(\int_{2^{k+1}B}|(b_{2})_{2B}-b_2(y_2)||f_{2}(y_2)|d\mu(y_2)\bigg)d\mu(y_1)\\
&&+\sum^{\infty}_{j=1}\int_{2^{j+1}B\backslash2^{j}B}\frac{|(b_{2})_{2B}-b_2(y_2)||f_{2}(y_2)|}{[V(x,y_{2})]^{2}}\\
&&\times \bigg(\int_{2^{j+1}B}|(b_1)_{2B}-b_1(y_1)||f_1(y_1)|d\mu(y_1)\bigg)d\mu(y_2)\\
&\lesssim& \|b_{2}\|_{BMO(\mu)}\|f_2\|_{\mathcal{M}^{p_2,\phi}(\omega)}\sum^{\infty}_{k=1}\frac{k}{\mu(2^{k}B)}\int_{2^{k+1}B}|(b_1)_{2B}-b_1(y_1)||f_1(y_1)|d\mu(y_1)\\
&&\times \bigg(\frac{\phi(\omega(2^{k+1}B))}{\omega(2^{k+1}B)}\bigg)^{\frac{1}{p_2}}\\
&&+C \|b_{1}\|_{BMO(\mu)}\|f_1\|_{\mathcal{M}^{p_1,\phi}(\omega)}\sum^{\infty}_{j=1}\frac{j}{\mu(2^{j}B)}\int_{2^{j+1}B}|(b_2)_{2B}-b_2(y_2)||f_2(y_2)|d\mu(y_2)\\
&&\times \bigg(\frac{\phi(\omega(2^{j+1}B))}{\omega(2^{j+1}B)}\bigg)^{\frac{1}{p_1}}\\
&\lesssim& \|b_{2}\|_{BMO(\mu)}\|f_2\|_{\mathcal{M}^{p_2,\phi}(\omega)}\sum^{\infty}_{k=1}\frac{k}{\mu(2^{k}B)}\bigg(\int_{2^{k+1}B}|f_1(y_1)|^{p_1}\omega(y_1)d\mu(y_1)\bigg)^{\frac{1}{p_1}}\\
&&\times\bigg(\int_{2^{k+1}B}|(b_1)_{2B}-b_1(y_1)|^{p^{\prime}_1}[\omega(y_1)]^{-\frac{p^{\prime}_1}{p_1}}\bigg)^{\frac{1}{p^{\prime}_1}}\bigg(\frac{\phi(\omega(2^{k+1}B))}{\omega(2^{k+1}B)}\bigg)^{\frac{1}{p_2}}\\
&&+C \|b_{1}\|_{BMO(\mu)}\|f_1\|_{\mathcal{M}^{p_1,\phi}(\omega)}\sum^{\infty}_{j=1}\frac{j}{\mu(2^{j}B)}\bigg(\int_{2^{j+1}B}|f_2(y_2)|^{p_2}\omega(y_2)d\mu(y_2)\bigg)^{\frac{1}{p_2}}\\
&&\times\bigg(\int_{2^{j+1}B}|(b_2)_{2B}-b_2(y_2)|^{p^{\prime}_2}[\omega(y_2)]^{-\frac{p^{\prime}_2}{p_2}}\bigg)^{\frac{1}{p^{\prime}_2}}\bigg(\frac{\phi(\omega(2^{j+1}B))}{\omega(2^{j+1}B)}\bigg)^{\frac{1}{p_1}}\\
&\lesssim& \|b_{1}\|_{BMO(\mu)}\|b_{2}\|_{BMO(\mu)}\|f_1\|_{\mathcal{M}^{p_1,\phi}(\omega)}\|f_2\|_{\mathcal{M}^{p_2,\phi}(\omega)}\bigg[\frac{\phi(\omega(B))}{\omega(B)}\bigg]^{\frac{1}{p}},
\end{eqnarray*}
which, combing the estimates of $J_{41}$, $J_{42}$ and $J_{43}$, Lemma 2.4, implies that
\begin{eqnarray*}
&&\|[b_1,b_2,T_{\theta}](f^{\infty}_1,f^{\infty}_2)\|_{\mathcal{M}^{p,\phi}(\omega)}\\
&&\quad=\mathop{\sup}\limits_{B}\frac{1}{\phi(\omega(B))^{\frac{1}{p}}}\bigg(\int_{B}|[b_1,b_2,T_{\theta}](f^{\infty}_1,f^{\infty}_2)(x)|^{p}\omega(x)d\mu(x)\bigg)^{\frac{1}{p}}\\
&&\quad\lesssim \|b_1\|_{BMO(\mu)}\|b_2\|_{BMO(\mu)}\|f_{1}\|_{\mathcal{M}^{p_1,\phi}(\omega)}\|f_2\|_{\mathcal{M}^{p_2,\phi}(\omega)}
\bigg[\frac{\omega(B)}{\phi(\omega(B))}\bigg]^{\frac{1}{p}}\bigg[\frac{\phi(\omega(B))}{\omega(B)}\bigg]^{\frac{1}{p}}\\
&&\quad\quad+\|b_{2}\|_{BMO(\mu)}\|f_{1}\|_{\mathcal{M}^{p_1,\phi}(\omega)}\|f_2\|_{\mathcal{M}^{p_2,\phi}(\omega)}\bigg[\frac{\phi(\omega(B))}{\omega(B)}\bigg]^{\frac{1}{p}}\\
&&\quad\quad\times\mathop{\sup}\limits_{B}\frac{1}{\phi(\omega(B))^{\frac{1}{p}}}\bigg(\int_{B}|b_1(x)-(b_1)_{2B}|^{p}\omega{(x)}d\mu(x)\bigg)^{\frac{1}{p}}\\
&&\quad\quad+\|b_{1}\|_{BMO(\mu)}\|f_{1}\|_{\mathcal{M}^{p_1,\phi}(\omega)}\|f_2\|_{\mathcal{M}^{p_2,\phi}(\omega)}\bigg[\frac{\phi(\omega(B))}{\omega(B)}\bigg]^{\frac{1}{p}}\\
&&\quad\quad\times\mathop{\sup}\limits_{B}\frac{1}{\phi(\omega(B))^{\frac{1}{p}}}\bigg(\int_{B}|b_2(x)-(b_2)_{2B}|^{p}\omega{(x)}d\mu(x)\bigg)^{\frac{1}{p}}\\
&&\quad\quad+\|f_{1}\|_{\mathcal{M}^{p_1,\phi}(\omega)}\|f_2\|_{\mathcal{M}^{p_2,\phi}(\omega)}\bigg[\frac{\phi(\omega(B))}{\omega(B)}\bigg]^{\frac{1}{p}}\\
&&\quad\quad\times\mathop{\sup}\limits_{B}\frac{1}{\phi(\omega(B))^{\frac{1}{p}}}\bigg(\int_{B}|b_1(x)-(b_1)_{2B}|^{p}|b_2(x)-(b_2)_{2B}|^{p}\omega{(x)}d\mu(x)\bigg)^{\frac{1}{p}}\\
&&\quad\lesssim \|b_1\|_{BMO(\mu)}\|b_2\|_{BMO(\mu)}\|f_{1}\|_{\mathcal{M}^{p_1,\phi}(\omega)}\|f_2\|_{\mathcal{M}^{p_2,\phi}(\omega)}\\
&&\quad\quad+\|f_{1}\|_{\mathcal{M}^{p_1,\phi}(\omega)}\|f_2\|_{\mathcal{M}^{p_2,\phi}(\omega)}\mathop{\sup}\limits_{B}\bigg[\frac{\phi(\omega(B))}{\omega(B)}\bigg]^{\frac{1}{p}}\frac{1}{\phi(\omega(B))^{\frac{1}{p}}}\\
&&\quad\quad\times\bigg(\int_{B}|b_1(x)-(b_1)_{2B}|^{p_{1}}\omega(x)d\mu(x)\bigg)^{\frac{1}{p_1}}\bigg(\int_{B}|b_2(x)-(b_2)_{2B}|^{p_{2}}\omega{(x)}d\mu(x)\bigg)^{\frac{1}{p_2}}\\
&&\quad\lesssim \|b_1\|_{BMO(\mu)}\|b_2\|_{BMO(\mu)}\|f_{1}\|_{\mathcal{M}^{p_1,\phi}(\omega)}\|f_2\|_{\mathcal{M}^{p_2,\phi}(\omega)}.
\end{eqnarray*}
Hence, the proof of the Theorem 1.9 is completed.

\noindent$\textbf{Conflict of interest}$\\
The authors declare that there is no conflict of interests.

\noindent$\textbf{Acknowledgments}$ \\
This research was supported by National Natural Science Foundation of China(Grant No. 11561062).



\begin{thebibliography}{999}

\bibitem{A} D. R. Adams, A note on Riesz potentials, Duke Math. J. {\bf42} (4) (1975), 765-778.

\bibitem{CW1} R. R. Coifman, G. Weiss, Extensions of Hardy spaces and their use in analysis, Bull. Am. Math. Soc. {\bf 83} (4) (1997), 569-645.

\bibitem{CW2} R. R. Coifman, G. Weiss, Analyse harmonique non-commutative sur certains espaces homog¨¨nes. Lecture Notes in Mathematics, {\bf 242}. Springer, Berlin (1971).

\bibitem{CL} J. Chou, X. Li, Y. Tong, H. Lin,  Generalized weighted Morrey spaces on RD-spaces, Rocky Mountain J. Math. {\bf 50} (4) (2020), 1277-1293.

\bibitem{CF} F. Chiarenza, M. Frasca, Morrey spaces and Hardy-Littlewood maximal function, Rend. Math. {\bf7} (1987), 273-279.

\bibitem{DR} J. Duoandikoetxea, M.Rosenthal, Extension and boundedness of operators on Morrey spaces from extrapolation techniques and embeddings. J. Geom. Anal. {\bf 28} (4) (2018), 3081-3108.

\bibitem{DXY} X. T. Duong, J. Xiao, L. Yan, Old and new Morrey spaces with heat kernel bounds. J. Fourier Anal. Appl. {\bf 13} (1) (2007), 87-111.

\bibitem{G} M. Giaquinta, Multiple integrals in the calculus of variations and nonlinear elliptic systems, Princeton Univ. Press, Princeton, NJ, 1983.

\bibitem{GT} D. Gilbarg, N. S. Trudinger, Elliptic Partial Differential Equations of Second Order,2nd ed., Springer-Verlag, Berlin, (1983).

\bibitem{GLMY} L. Grafakos, L. Liu, D. Maldonado, D. Yang, Multilinear analysis on metric spaces, Dissertationes Math. {\bf 497} (2014)

\bibitem{GLY} L. Grafakos, L. Liu, D. Yang, Boundedness of paraproduct operators on RD-spaces. Sci. China Math. {\bf 53} (8) (2010), 2097-2114.

\bibitem{HZ}
S. He, J. Zhou, Vector-valued maximal multilinear Calder\'on-Zygmund operator with nonsmooth kernel on weighted Morrey space. J. Pseudo-Differ. Oper. Appl.
{\bf 8} (2) (2017), 213-239.

\bibitem{JYY} X. Jiang, D. Yang, W. Yuan, Real interpolation for grand Besov and Triebel-Lizorkin spaces on RD-spaces. Ann. Acad. Sci. Fenn. Math. {\bf 36} (2) (2011), 509-529.

\bibitem{KS} Y. Komori and S. Shirai, Weighted Morrey spaces and a singular integral operator, Math. Nachr. {\bf 282} (2) (2009), 219-231.

\bibitem{P} V. Paatashvili, Boundary value problems for snalytic and harmonic functions in nonstandard spaces, Nova Science Publishers Inc., New York (2012)

\bibitem{KMRS} V. Kokilashvili, A. Meskhi, H. Rafeiro, S. Samko, Integral operators in non-standard function spaces: variable exponent H\"{o}lder, Morrey¨CCampanato and grand spaces, {\bf 2}. Birk$\ddot{a}$user/Springer,Heidelberg (2016)

\bibitem{LY} H. Lin, D. Yang, Pointwise multipliers for localized Morrey-Campanato spaces on RD-spaces. Acta Math. Sci. Ser. B (Engl. Ed.) {\bf 34} (6) (2014), 1677-1694.

\bibitem{LLW} Q. Li, H. Lin, X. Wang, Boundedness of commutators of $\theta$-type Calder\'on-Zygmund operators on generalized weighted Morrey spaces over RD-spaces, Anal. Math. Phys. {\bf 12} (1) (2022), 1-27.

\bibitem{M} C. B. Morrey, On the solutions of quasi-linear elliptic partial differential equations, Trans. Amer. Math. Soc. {\bf43} (1938),126-166.

\bibitem{M1} B. Muckenhoupt, Weighted norm in equalities for the Hardy maximal function. Trans. Am. Math. Soc. {\bf 165} (1972), 207-226.

\bibitem{N} E. Nakai, The Campanato, Morrey and H\"{o}lder spaces on spaces of homogeneous type. Studia Math. {\bf 176} (1) (2006), 1-19.

\bibitem{P1} J. Peetre,  On the theory of $L^{p,\lambda}$ spaces, J. Funct. Anal. {\bf 4}(1969), 71-87.

\bibitem{S} N. Samko, Weighted Hardy and singular operators in Morrey spaces. J. Math. Anal. Appl. {\bf 350} (1) (2009), 56-72.

\bibitem{SE} Y. Sawano, S. R. El-Shabrawy, Weak Morrey spaces with applications. Math. Nachr. {\bf 291} (1) (2018), 178-186.

\bibitem{ST} Y. Sawano, H. Tanaka, Morrey spaces for non-doubling measures. Acta Math. Sin. (Engl. Ser.) {\bf 21} (6) (2005), 1535-1544.

\bibitem{ST1} O. Str$\rm\ddot{o}$mberg, A. Torchinsky, Weighted Hardy Spaces, Lecture Notes in Mathematics, {\bf1381}. Springer, Berlin (1989).

\bibitem{W} H. Wang, Boundedness of $\theta$-type Calder\'on-Zygmund operators and commutators in the generalized weighted Morrey spaces. J. Funct. Spaces Art. ID 1309348 (2016)

\bibitem{YZ} D. Yang, Y. Zhou, Boundedness of sublinear operators in Hardy spaces on RD-spaces via atoms. J. Math. Anal. Appl. {\bf 339} (1) (2008), 622-635.

\bibitem{ZLL} S. Zhang, H. Lin, Y. Lin, The weighted Morrey boundedness of multilinear singular integral operators on RD-spaces. Anal. Theory Appl. {\bf 37} (3) (2021), 465-480.

\end{thebibliography}
\end{document}